\documentclass[12pt,leqno]{article}
\usepackage{amsmath, amsthm, amsfonts, amssymb, color}
\usepackage{mathrsfs}
\theoremstyle{definition}

\newcommand{\scr}[1]{\mathscr #1}
\definecolor{wco}{rgb}{0.5,0.2,0.3}

\numberwithin{equation}{section} \theoremstyle{remark}

\newcommand{\ua}{\uparrow}

\title{
{\bf Equivalent Harnack and Gradient Inequalities for Pointwise Curvature Lower Bound}
\footnote{Supported in part by NNSFC(11131003), SRFDP, the Laboratory of Mathematical and  Complex Systems and the Fundamental Research Funds for the Central Universities.}
}
\author{
{\bf Marc Arnaudon$^{3)}$, Anton Thalmaier$^{2)}$, Feng-Yu Wang$^{1),4)}$ }\\ \\
\footnotesize{$^1)$ School of Mathematical Sciences, Beijing Normal University, Beijing 100875, China}\\
\footnotesize{Email: \tttext{wangfy@bnu.edu.cn}}\\
 \footnotesize{$^2)$ Mathematics Research Unit, FSTC, University of Luxembourg}\\
 \footnotesize{6 rue Richard Coudenhove-Kalergi, L--1359 Luxembourg, Grand-Duchy of
  Luxembourg} \\
 \footnotesize{Email: \tttext{anton.thalmaier@uni.lu}}\\
  \footnotesize{$^3)$ Laboratoire de Math\'ematiques et Applications (UMR7348) Universit\'e de Poitiers,}\\
 \footnotesize{T\'el\'eport 2 - BP 30179 F--86962 Futuroscope Chasseneuil Cedex, France}\\
  \footnotesize{Email: \tttext{marc.arnaudon@math.univ-poitiers.fr}}\\
  \footnotesize{$^4)$ Department of Mathematics, Swansea University, Singleton Park, SA2 8PP, UK}\\
\footnotesize{Email: \tttext{F.Y.Wang@swansea.ac.uk}} }

\begin{document}
\def\tttext#1{{\normalfont\ttfamily#1}}
\def\R{\mathbb R}  \def\ff{\frac} \def\ss{\sqrt} \def\B{\mathbf B}
\def\N{\mathbb N} \def\kk{\kappa} \def\m{{\bf m}}
\def\dd{\delta} \def\DD{\Delta} \def\vv{\varepsilon} \def\rr{\rho}
\def\<{\langle} \def\>{\rangle} \def\GG{\Gamma} \def\gg{\gamma}
  \def\nn{\nabla} \def\pp{\partial} \def\EE{\scr E}
\def\d{\text{\rm{d}}} \def\bb{\beta} \def\aa{\alpha} \def\D{\scr D}
  \def\si{\sigma} \def\ess{\text{\rm{ess}}}
\def\beg{\begin} \def\beq{\begin{equation}}  \def\F{\scr F}
\def\Ric{\text{\rm{Ric}}} \def\Hess{\text{\rm{Hess}}}
\def\e{\text{\rm{e}}} \def\ua{\underline a} \def\OO{\Omega}  \def\oo{\omega}
 \def\tt{\tilde} \def\Ric{\text{\rm{Ric}}}
\def\cut{\text{\rm{cut}}} \def\P{\mathbb P}
\def\C{\scr C}     \def\E{\mathbb E}
\def\Z{\mathbb Z} \def\II{\mathbb I}
  \def\Q{\mathbb Q}  \def\LL{\Lambda}
  \def\B{\scr B}    \def\ll{\lambda}
\def\vp{\varphi}\def\H{\mathbb H}\def\ee{\mathbf e}

\maketitle
\begin{abstract} By using a coupling method, an explicit log-Harnack inequality with local geometry quantities  is established for (sub-Markovian) diffusion semigroups on a Riemannian manifold (possibly with boundary). This inequality as well as the consequent   $L^2$-gradient inequality, are proved to be equivalent to  the pointwise curvature lower bound condition together with the convexity or absence of the boundary.  Some applications of the log-Harnack inequality are also introduced.
\end{abstract} \noindent

 AMS subject Classification:\ 58J65, 60H30.   \\
\noindent
 Keywords: Log-Harnack inequality, Riemannian manifold, diffusion process.
 \vskip 2cm

\section{Introduction} Let $M$ be a $d$-dimensional connected complete Riemannian manifold possibly with a   boundary $\pp M$. Consider   $L=\DD+Z$ for a $C^1$-vector field~$Z$. Let $X_t(x)$ be the  (reflecting) diffusion process generated by $L$ with starting point $x$ and life time $\zeta(x)$. Then the associated diffusion semigroup $P_t$ is given by
 $$P_tf(x):= \E \big[f(X_t(x))1_{\{t<\zeta(x)\}}\big],\ \ t\ge 0, f\in \B_b(M).$$ Although the semigroup depends on $Z$ and the geometry on the whole manifold, we aim to establish  Harnack, resp.~gradient type  inequalities for $P_t$ by using local geometry quantities.

 Let $K\in C(M)$ be such that
 \beq\label{C} \Ric_Z:=\Ric-\nn Z\ge -K,\end{equation} i.e. for any $x\in M$ and $X\in T_xM$, $\Ric(X,X)- \<X,\nn_X Z\>\ge -K(x)|X|^2.$
 Next, for any   $D\subset M$, let
 $$K(D):= \sup_D K,\ \ D_r=\{z\in M: \rr(z,D)\le r\},\ \ r\ge 0,$$ where $\rr$ is the Riemannian distance on $M$.
 Finally, to investigate $P_t$ using   local curvature bounds, we introduce, for a given bounded   open domain $D\subset M$, the following class of reference functions:
 $$\C_D=\big\{\phi\in C^2(\bar D):\  \phi|_D>0,\ \phi|_{\pp D\setminus\pp M}=0,\ N\phi|_{\pp M\cap\pp D}\ge 0\big\},$$ where $N$ is the inward unit normal vector field of $\pp M$. When $\pp M=\emptyset,$ the restriction $N\phi|_{\pp M}\ge 0$ is automatically dropped.
 For any $\phi\in \C_D,$  we have
$$c_D(\phi)= \sup_D\big\{5|\nn \phi|^2  -\phi L\phi\big\} \in [0,\infty).$$ The finiteness of $c_D(\phi)$ is trivial since $\bar D$ is compact. To see that $c_D(\phi)\ge 0$, we consider the following two situations:\beg{enumerate}

\item[(a)] There exists $x\in \pp D\setminus \pp M$. We have $\phi(x)=0$  so that $c_D(\phi)\ge \big\{5|\nn \phi|^2  -\phi L\phi\big\}(x)=0.$
 \item[(b)]
   When $\pp D\setminus \pp M=\emptyset$, we have $\bar D=M$. Otherwise, there exists $z\in M\setminus (D\cup\pp M)$, For any $z'\in D\setminus \pp M$, let $\gg: [0,1]\to M\setminus\pp M$ be a smooth curve linking $z$ and $z'.$ Since $z'\in D$ but $z\notin D$, there exists $s\in [0,1]$ such that $\gg(s)\in\pp D$. This is however impossible since $\pp D\subset \pp M$ and $\gg(s)\notin \pp M.$
  Therefore, in this case $M=\bar D$ is compact so that the reflecting diffusion process is non-explosive. Now, let $x\in \bar D$ such that $\phi(x)=\max_{\bar D}\phi$.   Since $N\phi|_{\pp M}\ge 0$ due to $\phi\in \C_D$, $\phi(X_t)-\phi(x) -\int_0^t L \phi(X_s)\,\d s$ is a sub-martingale so that $$\phi(x)\ge \E\phi(X_t) \ge \phi(x)+\int_0^t \E L\phi(X_s)\,\d s,\ \ t\ge 0.$$ This implies $L\phi(x)\le 0$ (known as the maximum principle) and thus, $$c_D(\phi)\ge \big\{5|\nn \phi|^2  -\phi L\phi\big\}(x)\ge 0.$$\end{enumerate}

 \beg{thm}\label{T1.1} Let $K\in C(M)$. The following statements are equivalent:
\beg{enumerate} \item[$(1)$]  $(\ref{C})$  holds and $\pp M$ is either empty or convex.
 \item[$(2)$] For any   bounded open domain $D\subset M$ and any $\phi\in\C_D$,      the log-Harnack inequality
 \beg{equation*} \beg{split}&P_T\log f(y)- \log(P_T f(x)+1-P_T1(x))\\
 &\ \le   \ff {\rr(x,y)^2} 2 \bigg(\ff {K(D_{\rr(x,y})}{1-\e^{-2K(D_{\rr(x,y)})T}} + \ff{ c_D(\phi)^2(\e^{2K (D_{\rr(x,y)})T}-1) }{ 2K(D_{\rr(x,y)})\phi(y)^4}\bigg),\  \ T>0,\ y\in D,\  x\in M,   \end{split}\end{equation*} holds for strictly positive   $f\in \B_b(M).$
 \item[$(3)$]  For any   bounded open domain $D\subset M$ and any $\phi\in\C_D$,
 $$|\nn P_T f|^2(x)\le \big\{P_Tf^2-(P_Tf)^2\big\} (x) \bigg(\ff{K(D)}{1-\e^{-2K(D)T}} +\ff{ c_D(\phi)^2(\e^{2K (D)T}-1) }{ 2K(D)\phi(x)^4}\bigg)$$ holds for all $ x\in D, T>0,  f\in \B_b(M).$   \end{enumerate} If moreover $P_T1=1$, then the statements above are also equivalent to \beg{enumerate}
 \item[$(4)$]  For any   bounded open domain $D\subset M$ and any $\phi\in\C_D$, the Harnack type inequality
$$P_T f(y)\le P_Tf(x) + \rr(x,y)  \bigg(\ff {K(D)}{1-\e^{-2K(D)T}} + \ff{ c_D(\phi)^2(\e^{2K (D)T}-1) }{ 2K(D)\inf_{\ell(x,y)}\phi^4}\bigg)^{1/2}\ss{P_Tf^2(y)}$$ holds for  nonnegative $f\in \B_b(M)$,  $T>0$ and $x,y\in D$ such that the minimal geodesic $\ell(x,y)$ linking $x$ and $y$ is contained in $D$.
 \end{enumerate} \end{thm}
\paragraph{Remark} (i) When $K$ is constant, a number of equivalent semigroup inequalities are available for the curvature condition (\ref{C}) together with the convexity or absence of the boundary, see \cite{W10, WDM} and references within (see also \cite{BL,W} for equivalent semigroup inequalities of the curvature-dimension condition).  When $\pp M$ is either empty or convex, the above result provides at the first time equivalent semigroup properties for the general pointwise curvature lower bound condition.

(ii) When the diffusion process is explosive, the appearance of $1-P_T1$ in the log-Harnack inequality is essential.  Indeed, without this term the inequality does not hold for e.g. $f\equiv 1$ provided $P_T 1<1.$

(iii) The following result shows that the constant $1/2$ involved in the log-Harnack inequality  is sharp.

\beg{prp} \label{P1.2} Let $c>0$ be a constant. For any $x\in M$, strictly positive function $f $ with $|\nn f|(x)>0$ and $\log f\in C_0^2(M)$, and any constants $C>0 $, the   inequality
$$P_T \log f(y)\le \log\left(P_T f(x)+1-P_T1(x)\right) +  c\,\rr(x,y)^2 \bigg( \ff{C}{1-\e^{-2CT}}  + {\rm o}\Big(\ff 1 T\Big)\bigg)$$ for  small $T>0$ and  small  $\rr(x,y)$   implies that $c\ge 1/2.$
\end{prp}

\beg{proof} Let us take $v\in T_x M$ and $y_s=\exp_x[s v], s\ge 0.$ Then the given log-Harnack inequality implies that
\beq\label{K1}P_s \log f(y_s)- \log \left(P_s f(x)+1-P_s1(x)\right)\le cs^2|v|^2 \bigg( \ff{C}{1-\e^{-2Cs}}  + {\rm o}\Big(\ff 1 s\Big)\bigg) \end{equation}  holds for small $s>0.$ On the other hand,  for any $g\in C^2(M)$ with bounded $Lg$,  one has
\beq\label{K}\ff{\d}{\d s}P_s g|_{s=0}= Lg.\end{equation} Indeed, letting $X_t$ be the   diffusion process generated by $L$ with $X_0=x$,  by   It\^o's formula and the  dominated convergence theorem we obtain
$$\lim_{s\downarrow 0}\ff{P_s g(x)-g(x)}s = \lim_{s\downarrow  0} \ff 1 s \E\int_0^{s\land\zeta(x)}Lg(X_r)\,\d r =\E \lim_{s\downarrow  0} \ff 1 s \int_0^{s\land\zeta(x)} L g(X_r)\,\d r= Lg(x).$$    Combining (\ref{K1}) with (\ref{K}) we obtain
\beg{equation*}\beg{split}  \<v,&\nn  \log f\>(x)- |\nn\log f|^2(x)= L\log f(x) + \<v, \nn \log f\>(x) -\ff{Lf(x)}{f(x)} \\
&=\lim_{s\downarrow 0} \ff 1 s \big\{P_s \log f(y_s)-\log (P_s f(x)+1-P_s 1(x))\big\} \le \ff{c|v|^2} 2.\end{split}\end{equation*}
Taking $v= r\nn\log f(x)$ for $r\ge 0$ we obtain
$$\Big(r-1-\ff{cr^2}2\Big)|\nn\log f(x)|^2 \le 0,\ \ r\ge 0.$$ This implies $c\ge 1/ 2$ by taking $r=1/ c.$
\end{proof}

To derive the explicit log-Harnack inequality using local geometry quantities, we may take e.g. $D=B(y,1):= \{z: \rr(y,z)<1\}$. Let
\beg{equation*}\beg{split} K_y&=0\lor K(B(y,1)),\ \  K_{x,y} =  K(B(y,1+\rr(x,y))),\\
K_y^0&= 0\lor \sup\big\{-\Ric(U,U): U\in T_zM,\ |U|=1,\  z\in B(y,1)\big\},\\
b_y&= \sup_{B(y,1)} |Z|.\end{split}\end{equation*} Then $K(D_{\rr(x,y)})= K_{x,y}$ and  according to \cite[Proof of Corollary 5.1]{TW} (see page 121 therein with $\bar\dd_x$ replaced by 1), we may take $\phi(z) =\cos\ff{\pi \rr(y,z)}{2}$ so that $\phi(y)=1$ and
$$\kk(y):= K_y +\ff{\pi^2(d+3)}4 + \pi\Big(b_y+\ff 1 2  \ss{K_y^0(d-1)}\Big)\ge c_D(\phi).$$ Note that when $\pp M$ is convex, $N\rr(\cdot, y)|_{\pp M}\le 0$ so that $N\phi|_{\pp D\cap\pp M}\ge 0$ as required in the definition of $\C_D$.
Therefore,   Theorem \ref{T1.1} (2) implies  that
 \beg{equation}\label{LH2} P_t\log f(y)\le \log\big\{P_t f(x)+1-P_t1(x)\big\}+
   \ff{\rr(x,y)^2}2 \bigg(\ff{ K_{x,y}}{ 1-\e^{-2K_{x,y}t}} +\ff{ \kk(y)^2(\e^{2K_{x,y}t}-1)}{2K_{x,y}} \bigg) \end{equation}
 holds for all strictly positive $f\in \B_b(M), x,y\in M$ and $t>0.$ As in the proofs of  \cite[Corollary 1.2]{RW10} and \cite[Corollary 1.3]{W11}, this implies the following heat kernel estimates and entropy-cost inequality. When $P_t$ obeys the log-Sobolev inequality for $t>0$, the second inequality in Corollary \ref{C1.3}(2) below also implies the HWI inequality as shown in \cite{BGL, OV}. 

 \beg{cor} \label{C1.3} Assume $(\ref{C})$ and that $\pp M$ is either convex or empty. Let $Z=\nn V$ for some $V\in C^2(M)$ such that  $P_t$ is  symmetric w.r.t. $\mu(\d x):= \e^{V(x)}\,\d x$, where $\d x$ is the volume measure. Let $p_t$ be the density of $P_t$ w.r.t. $\mu$. Assume that   $(\ref{C})$ holds.
 \beg{enumerate} \item[$(1)$] Let $\bar K(y)= K(B(y,2)).$ Then
\beg{equation*}\beg{split} &\int_M p_t(y,z)\log p_t(y,z)\mu(\d z)\\
&\le \ss{t\land 1} \bigg(\ff{\bar K(y)}{1-\e^{-2\bar K(y)t}} +\ff{ \kk(y)^2(\e^{2\bar K(y)t}-1)}{2\bar K(y)} \bigg)+\log \ff{P_{2t}1(y)+\mu(1-P_t 1)}{\mu(B(y,\ss{t\land 1}))}\end{split}\end{equation*} holds for all $y\in M$ and $t>0.$
 \item[$(2)$] If $\mu$ is a probability measure and $P_t 1=1$, then the Gaussian heat kernel lower bound
$$ p_{2t} (x,y)\ge \exp\bigg[- \ff{\rr(x,y)^2}2 \bigg(\ff{ K_{x,y}}{ 1-\e^{-2K_{x,y}t}} +\ff{ \kk(y)^2(\e^{2K_{x,y}t}-1)}{2K_{x,y}} \bigg)\bigg],\ \ t>0,\  x,y\in M,$$ and the entropy-cost inequality
\beg{equation*}\beg{split} & \int_M (P_t f)\log P_t f\,\d\mu\\
& \le \inf_{\pi\in \C(\mu,f\mu)}\int_{M\times M} \ff{\rr(x,y)^2}2 \bigg(\ff{ K_{x,y}}{ 1-\e^{-2K_{x,y}t}} +\ff{ \kk(y)^2(\e^{2K_{x,y}t}-1)}{2K_{x,y}} \bigg)\pi(\d x,\d y),\ \ t>0,\end{split}\end{equation*}  hold for any probability density function  $f$ of $\mu$, where $\C(\mu,f\mu)$ is the set of all couplings of $\mu$ and $f\mu$. \end{enumerate}\end{cor}
\beg{proof} According to (\ref{LH2}), the heat kernel lower bound in (2) follows from the proof of \cite[Corollary 1.3]{W11}, while the other two inequalities can be proved as in the proof of \cite[Corollary 1.2]{RW10}. Below we only present a brief proof of (1).

By an approximation argument we may apply (\ref{LH2}) to $f(z):= p_t(y,z)$ so that
\beg{equation*}\beg{split} I& := \int_M p_t(y,z)\log p_t(y,z)\mu(\d z)\\
&\le \log \{p_{2t}(x,y)+1-P_t 1(x)\} + \ff{\rr(x,y)^2}2 \bigg(\ff{ K_{x,y}}{ 1-\e^{-2K_{x,y}t}} +\ff{ \kk(y)^2(\e^{2K_{x,y}t}-1)}{2K_{x,y}} \bigg) .\end{split}\end{equation*} Since $K_{x,y}\le \bar K(y)$ for $x\in B(y,1)$, this implies
that
\beg{equation*}\beg{split} &\e^I \mu\Big(B\big(y,\ss{t\land 1}\big)\Big)\exp\bigg[-\ff {t\land 1}2   \bigg(\ff{ \bar K(y)}{ 1-\e^{-2\bar K(y)t}} +\ff{ \kk(y)^2(\e^{2\bar K(y)t}-1)}{2\bar K(y)} \bigg) \bigg]\\
&\le \e^I \int_M \exp\bigg[- \ff{\rr(x,y)^2}2 \bigg(\ff{ K_{x,y}}{ 1-\e^{-2K_{x,y}t}} +\ff{ \kk(y)^2(\e^{2K_{x,y}t}-1)}{2K_{x,y}} \bigg) \bigg]\mu(\d x)\\
&\le \int_M \{p_{2t}(x,y)+1-P_t1(x)\}\mu(\d x)=P_{2t} 1(y)+\mu(1-P_t 1).\end{split}\end{equation*} This proves (1). \end{proof}
We remark that    the entropy upper bound in (1) is sharp for short time, since both $-\log \mu(B(y,\ss t))$ and  the entropy of the Gaussian heat kernel behave  like  $\ff d 2 \log \ff 1 t$ for small $t>0.$

\section{Proof of Theorem \ref{T1.1}}

 We first observe that when $P_T1=1$  the equivalence of (3) and (4) is implied by the proof of \cite[Proposition 1.3]{WN}.  Indeed, by (3)
  $$|\nn P_T f|^2\le \big\{P_Tf^2-(P_Tf)^2\big\}  \bigg(\ff{K(D)}{1-\e^{-2K(D)T}} +\ff{ c_D(\phi)^2(\e^{2K (D)T}-1) }{ 2K(D)\inf_{\ell(x,y)}\phi^4}\bigg)$$ holds on the minimal geodesic $\ell(x,y)$, so that the Harnack inequality in (4) follows from the first part in the proof of \cite[Proposition 1.3]{WN}. On the other hand, by the second part of the proof, the inequality in (4) implies
   $$|\nn P_t f|^2\le \big\{P_Tf^2\big\}  \bigg(\ff{K(D)}{1-\e^{-2K(D)T}} +\ff{ c_D(\phi)^2(\e^{2K (D)T}-1) }{ 2K(D) \phi^4}\bigg)$$ on $D$. Replacing $f$
   by $f-P_Tf(x)$, we obtain the inequality in (3) since $\nn P_T f= \nn P_T(f -P_T f(x))$ provided $P_T1=1.$

In the following three subsections, we prove (1) implying (2), (2) implying (3), and (3) implying (1) respectively.

\subsection{Proof of (1) implying (2)} We assume the curvature condition (\ref{C}) and that $\pp M$ is either empty or convex. To prove the log-Harnack inequality in~(2), we will make use of the coupling argument proposed in \cite{ATW06}.
 As explained in \cite[Section~3]{ATW06}, we may and do assume that the cut-locus of the manifold is empty.  

Now, let $T>0$ and $y\in D, x\ne y$ be fixed. For any $z,z'\in M$, let $P_{z,z'}\colon T_zM\to T_{z'}M$ be the parallel transport along the unique minimal geodesic from $z$ to $z'$.  Let $X_t$ solve the following It\^o type SDE on $M$
 $$\d^I X_t = \ss 2 \Phi_t \,\d B_t +Z(X_t) \,\d t+N(X_t)\d l_t,\ \ X_0=x,$$ up to the life time $\zeta(x)$, where $B_t$ is the $d$-dimensional Brownian motion,   $\Phi_t$ is the horizontal lift of $X_t$ on the frame bundle $O(M)$, and $\l_t$ is the local time of $X_t$ on $\pp M$ if $\pp M\ne\emptyset$.  When $\pp M=\emptyset$, we simply take $l_t=0$ so that the last term in the equation disappears. 
 
 To construct another process starting at $y$ such that it meets $X_t$ before $T$ and its hitting time to $\pp D$, let $Y_t$ solve the SDE with $ Y_0=y$
\beq\label{Y}\d^I Y_t = \sqrt2P_{X_t, Y_t}\Phi_t \,\d B_t+Z(Y_t)\,\d t - \ss{\xi_1(t)^2+\xi_2(t)^2}\,  \nn\rr(X_t,\cdot)(Y_t)\,\d t +N(Y_t)\d\tt l_t,\end{equation} where $\tt l_t$ is the local time of $Y_t$ on $\pp M$ when $\pp M\ne\emptyset$, and 
 \beg{equation*}\beg{split} &\xi_1(t)= \ff{2K(D_{\rr(x,y)})\exp[-K(D_{\rr(x,y)}) t]}{1-\exp[-2K(D_{\rr(x,y)}) T]}\rr(x,y)1_{\{Y_t\ne X_t\}},\\
 &\xi_2(t)= \ff{2c_D(\phi) \rr(X_t,Y_t)}{\phi(Y_t)^2}, \ \ \ t\in [0,T].\end{split}\end{equation*}
 Then $Y_t$ is well-defined before $T\land \tau_{D(x,y)}(x)\land \tau_D(y)$, where
 $$\tau_D(y):=\inf\{t\in [0, T\land \zeta(x)): Y_t\in \pp D\},\ \
 \tau_{D(x,y)}(x)= \inf\{t\ge 0: X_t\notin D(x,y)\}.$$  Let
 $$ \tau=\inf\{t\in [0,\zeta(x)\land\zeta(y)): X_t=Y_t\},$$ where $\inf\emptyset=\infty$ by convention.

Let $\Theta=  \tau\land T\land\tau_D(y)\land\tau_{D(x,y)}(x)$ and set
$$\eta(t)= \ff1{\sqrt2} \ss{\xi_1(t)^2+\xi_2(t)^2}\, \nn\rr(\cdot, Y_t)(X_t),\ \ t\in [0,\Theta).$$ Define
$$R= \exp\bigg[-\int_0^{\Theta}\<\eta(t),\Phi_t\,\d B_t\>-\ff 1 2\int_0^\Theta |\eta(t)|^2\,\d t\bigg].$$
We intend to prove
\beg{enumerate}\item[(i)] $R$ is a well-defined probability density   with
$$\E \big\{R\log R\big\} \le   \ff{\rr(x,y)^2}{2}\bigg(  \ff {K(D_{\rr(x,y})}{1-\e^{-2K(D_{\rr(x,y)})T}} + \ff{ c_D(\phi)(\e^{2K (D_{\rr(x,y)})^2T}-1) }{ 2K(D_{\rr(x,y)})\phi(y)^4}\bigg).$$
\item[(ii)] $\tau\le T\land\tau_D(y)\land\tau_{D(x,y)}(x)$ holds  $\Q$-a.s., where $\Q:= R\P.$  \end{enumerate}
Once these two assertions are confirmed, by taking $Y_t=X_t$ for $t\ge \tau$ we see that $Y_t$ solves (\ref{Y}) up to its life time $\zeta(y)=\zeta(x)$ and $X_T=Y_T$ for $T<\zeta(x)$. Moreover, by the Girsanov theorem the process
$$\tt B_t:= B_t + \int_0^t \eta(s)\,\d s,\ \ t\ge 0$$ is a $d$-dimensional Brownian motion under $\Q$ and equation (\ref{Y}) can be reformulated as
\beq\label{YY}\d^I Y_t = \ss 2 P_{X_t, Y_t}\Phi_t \,\d \tt B_t+Z(Y_t)\,\d t +N(Y_t) \d\tt l_t,\ \ Y_0=y.\end{equation} Combining this with the Young inequality
(see \cite[Lemma 2.4]{ATW09})
\beg{equation*} \beg{split} P_T \log f(y)&= \E\big\{R1_{\{T<\zeta(y)\}}\log f(Y_T)\big\}= \E\big\{R1_{\{T<\zeta(x)\}}\log f(X_T)\big\}\\
& \le \E R\log R + \log \E \exp[1_{\{T<\zeta(x)\}}\log f(X_T)]\\
&=\log(P_T f(x)+1-P_T1(x)) + \E R\log R\\
& \le \log(P_T f(x)+1-P_T1(x))\\
&\quad  +  \ff{\rr(x,y)^2}2\bigg(  \ff {K(D_{\rr(x,y})}{1-\e^{-2K(D_{\rr(x,y)})T}-1} + \ff{ c_D(\phi)^2(\e^{2K (D_{\rr(x,y)})T}-1) }{ 2K(D_{\rr(x,y)})\phi(y)^4}\bigg).\end{split}\end{equation*}
This gives the desired log-Harnack inequality.

Below we prove (i) and (ii) respectively.

\beg{lem}\label{L1} For any $n\ge 1,$ let $$\tau_n(y)= \inf\big\{t\in [0, T\land \zeta(x)):\ \rr(Y_t, D^c)\le n^{-1}\big\}$$ and $$\Theta_n=\tau\wedge \ff{nT}{n+1}\land \tau_{D(x,y)}(x)\land \tau_n(y).$$ Let $R_n$ be defined as $R$ using $\Theta_n$ in place of $\Theta$. Then $\{R_n\}_{n\ge 1}$ is a uniformly integrable martingale with $\E R_n=1$ and
$$ \E \{R_n\log R_n\} \le    \ff{\rr(x,y)^2}2\bigg(  \ff {K(D_{\rr(x,y})}{1-\e^{-2K(D_{\rr(x,y)})T}-1} + \ff{ c_D(\phi)^2(\e^{2K (D_{\rr(x,y)})T}-1) }{ 2K(D_{\rr(x,y)})\phi(y)^4}\bigg) $$ for $n\ge 1.$ Consequently, {\rm (i)} holds.
\end{lem}
\beg{proof} (i) follows from the first assertion and the martingale convergence theorem. Since before time $\Theta_n$ the process $\eta(t)$ is bounded, the martingale property and $\E R_n=1$ is well-known. So, it remains to prove the entropy upper bound.
By the It\^o formula we see that (cf. (2.3) and (2.4) in \cite{ATW06})
\beq\label{D2} \d \rr(X_t,Y_t)\le K(D_{\rr(x,y)}) \rr(X_t,Y_t)\,\d t -\ss{\xi_1(t)^2 +\xi_2(t)^2}\,\d t,\ \ t\le \Theta_n.\end{equation}
Then
$$\d \rr(X_t,Y_t)^2 \le  2K(D_{\rr(x,y)}) \rr(X_t,Y_t)^2\,\d t -\ff{4c_D(\phi) \rr(X_t,Y_t)^2}{\phi(Y_t)^2}\,\d t, \ \ \ t\le \Theta_n.$$
Note  that $(\tt B_t)_{t\in [0,\Theta_n]}$ is a $d$-dimensional Brownian motion under the probability $\Q_n:= R_n\P$. Combining this with (\ref{YY}) and using It\^o's formula along with the facts that  the martingale part of $\rr(X_t,Y_t)^2$ is zero and $N\phi|_{\pp D\cap\pp M}\ge 0$,  we obtain  \beg{equation*}\beg{split} \d \Big\{\ff {\rr(X_t,Y_t)^2}{\phi(Y_t)^4}\Big\}
&\le \d M_t  -\ff{4\rr(X_t,Y_t)^2}{\phi(Y_t)^6} \big\{c_D(\phi)+\phi(Y_t)L\phi(Y_t)-5|\nn\phi(Y_t)|^2\big\} \,\d t\\
&\qquad\qquad -\ff{2K(D_{\rr(x,y)})\rr(X_t,Y_t)^2}{\phi(Y_t)^4} \,\d t\\
& \le \d M_t -  \ff{2K(D_{\rr(x,y)})\rr(X_t,Y_t)^2}{\phi(Y_t)^4} \,\d t,
\  \   \ t\le \Theta_n,\end{split}\end{equation*} where
$$\d M_t:= -\ff{4\rr(X_t,Y_t)^2}{\phi(Y_t)^5}\<\nn \phi(Y_t), P_{X_t,Y_t}\Phi_t\,\d \tt B_t\>$$ is a $\Q_n$-martingale for $t\le \Theta_n$. This implies
$$\E_{\Q_n} \Big\{\ff{\rr(X_{t\land \Theta_n} ,Y_{t\land\Theta_n})^2}{\phi(Y_{t\land\Theta_n})^4}\Big\} \le \ff{\rr(x,y)^2}{\phi(y)^4}\,\e^{2K(D_{\rr(x,y)})t},\ \ t\ge 0.$$
 Hence,
\beg{equation*}\beg{split} &\E \big\{R_n\log R_n\big\} = \ff 1 2 \E_{\Q_n} \int_0^{\Theta_n} |\eta(t)|^2 \,\d t=\ff 1 4\E_{\Q_n}\int_0^{\Theta_n} \big\{\xi_1(t)^2+\xi_2(t)^2\big\}\,\d t\\
&\le \ff{ K(D_{\rr(x,y)})^2\rr(x,y)^2}{(1-\e^{-2K(D_{\rr(x,y)})T})^2}\int_0^T\e^{-2K(D_{\rr(x,y)})t}\,\d t+ c_D(\phi)^2\int_0^T \E_{\Q_n}   \ff{\rr(X_{t\land \Theta_n} ,Y_{t\land\Theta_n})^2}{\phi(Y_{t\land\Theta_n})^4} \,\d t\\
&\le \ff{  K(D_{\rr(x,y)}) \rr(x,y)^2}{2(1-\e^{-2K(D_{\rr(x,y)})T})} + \ff{ c_D(\phi)^2(\e^{2K( D_{\rr(x,y)})T} -1)\rr(x,y)^2}{2K(D_{\rr(x,y)})\phi(y)^4},\ \ s>0.\end{split}\end{equation*} \end{proof}

\beg{lem}\label{L2} We have $\tau\le T\land \tau_D(y)\land \tau_{D(x,y)}(x),\ \Q$-a.s.\end{lem}
\beg{proof} By (\ref{D2}) we have
\beq\label{DD} \int_0^{\Theta} \big\{\xi_1(t) + \xi_2(t)\big\}\,\d t=\lim_{n\to\infty} \int_0^{\Theta_n} \big\{\xi_1(t) + \xi_2(t)\big\}\,\d t<\infty.\end{equation}
Since under $\Q$ the process $Y_t$ is generated by $L$,  as observed in the beginning  of \cite[Section~4]{TW}  we have
$$\int_0^{\tau_D(y)} \ff{1}{\Phi(Y_t)^2}\,\d t=\infty,\ \ \Q\text{-a.s.}$$ Then (\ref{DD}) implies that $\Q$-a.s.
\beq\label{D3}\tau_D(y)> \tau_{D(x,y)}(x)\land\tau\land T.\end{equation} Moreover, it follows from (\ref{D2}) that
\beg{equation*}\beg{split} \rr(X_t,Y_t)&\le  \e^{K(D_{\rr(x,y)})t} \rr(x,y)-\int_0^t \e^{K(D_{\rr(x,y)})(t-s)}\xi_1(s)\,\d s \\
&\le \ff{\e^{-2K(D_{\rr(x,y)})t}- \e^{-2K(D_{\rr(x,y)})T}}{1-\e^{-2K(D_{\rr(x,y)})T}} \e^{K(D_{\rr(x,y)})t} \rr(x,y)\le \rr(x,y)1_{[0,T]}(t),\ \ t\in [0, \Theta_n].\end{split}\end{equation*}So, $\tau_{D(x,y)}\ge \tau_D(y) $ and $T\ge \tau.$ Combining these inequalities with (\ref{D3}) we complete the proof.
\end{proof}
\subsection{Proof of (2) implying (3)} We will present below a more general result, which works for sub-Markovian operators on metric  spaces.
Let $(E,\rr)$ be  a metric space, and let $P$ be a sub-Markovian operator on~$\B_b(E)$.
$$\dd(f)(x)=  \limsup_{y\to x} \ff{f(y)-f(x)}  {\rr(x,y)},\ \ x\in E, f\in \B_b(E).$$ If in particular $E=M$ and $f$ is differentiable at point $x$, then $\dd(f)(x)=|\nn f|(x).$  So, (2) implying (3) is a direct consequence of the following result.

\beg{prp} Let $x\in E$ be   fixed. If there exists a positive continuous function $\Phi$ on $E$ such that the log-Harnack inequality
\beq\label{LH0} P\log f(y)\le \log \big\{Pf(x) +1-P1(x)\big\}+\Phi(y) \rr(x,y)^2,\ \ f>0,\  f\in\B_b(E),\end{equation} holds for small $\rr(x,y)$, then
\beq\label{G0} \dd(Pf)^2(x)\le 2\Phi(x) \big\{Pf^2(x)-(Pf)^2(x)\big\},\ \ f\in\B_b(E).\end{equation} \end{prp}

\beg{proof} Let $f\in \B_b(E)$. According to the proof of \cite[Proposition 2.3]{W10}, (\ref{LH0}) for small $\rr(x,y)$ implies that $Pf$ is continuous at $x$.
Let $\{x_n\}_{n\ge 1}$ be a sequence converging to $x$, and denote $\vv_n= \rr(x_n,x)$. For any positive constant $c>0$, we apply (\ref{LH0}) to $c\vv_n f +1$ in place of $f$, so that for large enough $n$
$$P\log (c\vv_n f+1)(x_n)\le \log \big\{P(c\vv_n f +1)(x)+1 -P1 (x)\big\} + \Phi(x_n) \vv_n^2.$$ Noting that for large $n$ (or for small $\vv_n$) we have
\beg{equation*}\beg{split} &P\log (c\vv_n f+1)(x_n) = P\Big(c\vv_nf - \ff 1 2 (c\vv_n)^2 f^2\Big)(x_n)+ \text{o}(\vv_n^2)\\
&\qquad\qquad=c\vv_n Pf(x) + c\vv_n^2 \,\ff{Pf(x_n)-Pf(x)}{\rr(x_n,x)} -\ff 1 2 (c\vv_n)^2 Pf^2(x) +\text{o}(\vv_n^2),\\
& \log \big\{P(c\vv_n f +1)(x)+1 -P1 (x)\big\} = c\vv_n Pf(x) -\ff 1 2 (c\vv_n)^2 (Pf)^2(x) +\text{o}(\vv_n^2). \end{split}\end{equation*}
We obtain
$$c\limsup_{n\to\infty} \ff{Pf(x_n)- Pf(x)}{\rr(x_n,x)}\le \ff {c^2}  2 \big\{Pf^2(x)-(Pf)^2(x)\big\} + \Phi(x),\ \ c>0.$$ Therefore,
$$\dd(Pf)(x) \le \ff{c}2 \big\{Pf^2(x)-(Pf)^2(x)\big\} + \ff{\Phi(x)}c,\ \ c>0.$$ This implies (\ref{G0}) by minimizing the upper bound in $c>0.$
\end{proof}

\subsection{Proof of (3) implying (1)}  The proof of $\Ric_Z\ge -K$ is more or less standard by using the Taylor expansions for small $T>0$.  Let $x\in M\setminus \pp M$ and $D= B(x,r)\subset M\setminus \pp M$ for small $r>0$ such that $\phi:= r^2-\rr(x,\cdot)^2\in \C_D$. It is easy to see that for $f\in C_0^\infty(M)$ and small $t>0$,
\beg{equation*}\beg{split} &|\nn P_t f|^2(x)= |\nn f|^2(x)+ 2 t \<\nn f, \nn Lf\> +\text{o}(t),\\
&\ff{K(D)}{1-\e^{-2K(D)t} } +\ff{c_D(\phi)^2(\e^{2K(D)t}-1)}{2K(D)\phi(x)^4} = \ff 1 {2t} +\ff {K(D)} 2 +\text{o}(1).\end{split}\end{equation*} Moreover
(see \cite[(3.6)]{WDM}),
$$P_t f^2(x)-(P_tf)^2(x)= 2t |\nn f|^2(x) +t^2 \big\{2\<\nn f, \nn Lf\>+L|\nn f|^2\big\}(x) +\text{o}(t).$$ Combining these with (\ref{G0}) we obtain
$$ \GG_2(f)(x):= \ff 1 2 L|\nn f|^2(x) -\<\nn f,\nn Lf\>(x) \ge -K(D)|\nn f|^2(x)=-\big(\sup_{B(x,r)} K\big)|\nn f|^2(x).$$  Letting  $r\downarrow 0$, we arrive at $\GG_2(f)(x)\ge -K(x)$ for $x\in M\setminus \pp M$ and $f\in C_0^\infty(M)$, which is equivalent to  (\ref{C}).

Next, we assume that $\pp M\ne\emptyset$ and intend to prove from (3)  that the second fundamental form $\II$ of $\pp M$ is non-negative, i.e. $\pp M$ is convex.
When $M$ is compact, the proof was done in \cite{WDM} (see the proof of Theorem 1.1 therein for (7) implying  (1)). Below we show that the proof works for general setting by using a localization argument with a stopping time. 

Let $x\in \pp M$ and $r>0$. Define
$$\si_r=\inf\{s\ge 0:\rr(X_s,x)\ge r\},$$ where $X_s$ is the $L$-reflecting  diffusion process starting at point $x$. Let $l_s$ be the local time of the process on $\pp M$. Then, according to \cite[Lemmas 2.3 and 3.1]{W12}, there exist two constants $C_1,C_2>0$ such that
\beq\label{C1} \P(\si_r\le t)\le \e^{-C_1/t},\ \ t\in (0,1],\end{equation} and
\beq\label{C2}    \Big|\E l_{t\land\si_r}-\ff{2\ss t}{\ss\pi}\Big|\le C_2 t,\ \ t\in [0,1],\end{equation} where (\ref{C1}) is also ensured by \cite[Lemma 2.3]{ATW09} for $\pp M=\emptyset$. 
Let $f\in C_0^\infty(M)$ satisfy the Neumann boundary condition. We aim to prove $\II(\nn f,\nn f)(x)\ge 0.$ To apply Theorem \ref{T1.1}(3), we construct $D$ and $\phi\in \C_D$ as follows.

 Firstly, let $\varphi\in C_0^\infty(\pp M)$ such that $\varphi(x)=1$ and supp$\varphi\subset \pp M\cap B(x, r/2)$, where $B(x,s)=\{z\in M: \rr(z,x)<s\}$ for $s>0.$  Then
letting $\phi_0(\exp_y[sN])= \varphi(y)$ (where $y\in \pp M$, $s\ge 0$), we extend $\varphi$ to a smooth function in a neighborhood of $\pp M$, say $\pp_{r_0}M:= \{ z\in M: \rr(z,\pp M) <r_0\}$ for some $r_0\in (0,r)$  such that  $\rr(\cdot, \pp M)$ is smooth on $(\pp_{r_0}M)\cap B(x, r).$ Obviously, $\phi_0$ satisfies the Neumann boundary condition. Finally, for   $h\in C^\infty([0,\infty))$ with $h|_{[0,r_0/4]}=1$ and $h|_{[r_0/2,\infty)}=0$,  we take $\phi = \phi_0 h(\rr(\cdot,\pp M))$ and $D= \{z\in M: \phi(z)>0\}.$  Then $\phi(x)=1$, $\phi|_{\pp D\setminus \pp M}=0$, $N\phi|_{\pp M}=N\phi_0|_{\pp M}=0,$  and $D\subset B(x,r).$

Once $D$ and $\phi\in \C_D$ are given, below we calculate both sides of the gradient inequality in (3) respectively.

According to (\ref{C1}), for small $t>0$ we have
\beg{equation}\label{AA1}\beg{split} P_t f^2(x)&= \E f^2 (X_{t\land \si_r}) +\text{o}(t^2) =f^2(x) + \E\int_0^{t\land \si_r} L f^2(X_s)\,\d s + \text{o}(t^2)\\
&= f^2(x) + 2 \E\int_0^{t\land \si_r} (fL f)(X_s)\,\d s +  2\E\int_0^{t\land \si_r} |\nn f|^2(X_s) \,\d s + \text{o}(t^2).\end{split}\end{equation}
Noting that by the Neumann boundary condition
$$ \E|f(x)-f(X_{s\land\si_r})|^2 \le \|L(f(x)-f)^2\|_\infty s,\ \ s\ge 0,$$
we have
\beg{equation}\label{AA2}\beg{split} & \E\int_0^{t\land \si_r} (fL f)(X_s)\,\d s- f(x) \E\int_0^{t\land \si_r} L f(X_s)\,\d s\\
&=\E\int_0^{t\land \si_r} Lf(x)\{f(x)-f(X_{s})\}\,\d s
+ \E \int_0^{t\land\si_r} (Lf(X_s)-Lf(x)) (f(x)-f(X_s))\,\d s \\
& \le \|Lf\|_\infty \E \int_0^{t\land\si_r} \int_0^s  \,\d u  +\E\int_0^t \ss{\E|Lf(X_{s\land \si_r})-Lf(x)|^2 \cdot \E|f(x)-f(X_{s\land \si_r})|^2} \,\d s\\
&=\text{o}(t^{3/2}).\end{split}\end{equation}  Moreover, by the It\^o formula and the fact that $N|\nn f|^2 = 2 \II(\nn f,\nn f)$ holds on $\pp M$ (see e.g.~\cite[(3.8)]{WDM}), we have
\beg{equation*}\beg{split} \E |\nn f|^2(X_{s\land\si_r} ) &= |\nn f|^2(x) + \E\int_0^{s\land \si_r} L|\nn f|^2(X_u)\,\d u + 2 \int_0^{s\land \si_r} \II(\nn f,\nn f)(X_u) \,\d l_u\\
&\le  |\nn f|^2(x) + 2\II(r) \E l_{s\land \si_r} + \text{O}(t),\end{split}\end{equation*} where
$$\II(r):= \sup\big\{\II(\nn f,\nn f)(y):\ y\in \pp M\cap B(x,r)\big\}.$$ Combining this with (\ref{AA1}), (\ref{AA2}) and using (\ref{C1}) and (\ref{C2}), we obtain
\beg{equation}\label{C3} P_t f^2 (x) \le f^2(x)+ 2 f(x) \E\int_0^{t\land \si_r} L f (X_s)\,\d s + C t^{3/2} \II(r)+\text{o}(t^{3/2})\end{equation} for some constant $C>0$ and small $t>0.$

On the other hand, by (\ref{C1}) we have
$$(P_t f)^2(x)= \bigg(f(x)+ \E\int_0^{t\land\si_r} Lf(X_s)\,\d s +  \text{o}(t^2)\bigg)^2 = f^2(x)+ 2 t f(x)\E\int_0^{t\land\si_r} Lf(X_s)\,\d s + \text{o}(t^2).$$
Combining this with (\ref{C3}) and noting that
$$\ff{K(D)}{1-\e^{-2K(D)t}} +\ff{ c_D(\phi)^2(\e^{2K (D)t}-1) }{ 2K(D)\phi(x)^4}= \ff 1 {2 t} +\text{O}(1) $$ holds for small $t>0$, we arrive at
$$\big\{P_tf^2-(P_tf)^2\big\} (x) \bigg(\ff{K(D)}{1-\e^{-2K(D)t}} +\ff{ c_D(\phi)^2(\e^{2K (D)t}-1) }{ 2K(D)\phi(x)^4}\bigg)\le |\nn f|^2(x)+ C\II(r)\ss t  +\text{o}(t^{1/2})$$ for small $t>0.$ Combining this with the gradient inequality in (3) and noting that
$$|\nn P_t f|^2(x)= \bigg|\nn f(x)+ \int_0^t \nn P_s Lf(x)\,\d s\bigg|^2  =  |\nn f|^2(x) + \text{O}(t),$$ we conclude that
$$\II(r) \ge\lim_{t\to 0}  \ff 1{C\ss t} \bigg\{\big\{P_tf^2-(P_tf)^2\big\} (x) \bigg(\ff{K(D)}{1-\e^{-2K(D)t}} +\ff{ c_D(\phi)^2(\e^{2K (D)t}-1) }{ 2K(D)\phi(x)^4}\bigg)-|\nn P_t f|^2)(x)\bigg\}\ge 0.$$ Therefore,  $\II(\nn f,\nn f)(x)= \lim_{r\to 0}\II(r)\ge 0.$

\beg{thebibliography}{99}
 \bibitem{ATW06} M. Arnaudon, A. Thalmaier, F.-Y. Wang,
  \emph{Harnack inequality and heat kernel estimates
  on manifolds with curvature unbounded below,} Bull. Sci. Math. {\bf 130} (2006), 223--233.

\bibitem{ATW09} M. Arnaudon, A. Thalmaier, F.-Y. Wang,
  \emph{Gradient estimates and Harnack inequalities on non-compact Riemannian manifolds,}
   Stoch. Proc. Appl. {\bf 119} (2009), 3653--3670.
   
  \bibitem{BL} D. Bakry, M. Ledoux, \emph{A logarithmic Sobolev form of the Li-Yau parabolic inequality,}  Rev. Mat. Iberoam. {\bf 22} (2006),   683--702.
  
  \bibitem{BGL}  S. Bobkov, I. Gentil, M. Ledoux, \emph{Hypercontractivity of Hamilton-Jacobi
equations, }  J. Math. Pure Appl.    {\bf 80} (2001), 669--696.

  \bibitem{OV}  F. Otto, C. Villani,  \emph{Generalization of an inequality by Talagrand
and links with the logarithmic Sobolev inequality, }   J. Funct.
Anal.  {\bf 173}(2000), 361--400.

\bibitem{RW10} M. R\"ockner, F.-Y. Wang,  \emph{Log-Harnack  inequality for stochastic differential
  equations in Hilbert spaces and its consequences,}
  Inf. Dim. Anal. Quant. Proba. Relat. Top. {\bf 13} (2010), 27--37.

\bibitem{TW} A. Thalmaier, F.-Y. Wang, \emph{Gradient estimates for harmonic functions on regular domains in Riemannian manifolds,} J. Funct. Anal.  {\bf 155} (1998),  109--124.

\bibitem{W10} F.-Y. Wang, \emph{Harnack inequalities on manifolds with boundary and applications,}    J.
Math. Pures Appl.   {\bf 94} (2010), 304--321.

\bibitem{W11} F.-Y. Wang, \emph{Equivalent semigroup properties for the curvature-dimension condition,} Bull. Sci. Math. {\bf 135} (2011), 803--815.

\bibitem{WDM}   F.-Y. Wang, \emph{Semigroup properties for the second fundamental form,} Doc. Math. {\bf 15} (2010), 543--559.

\bibitem{W}   	F.-Y. Wang, \emph{Equivalent semigroup properties for the curvature-dimension condition,}  Bull. Sci. Math. {\bf 135} (2011) 803--815

\bibitem{WN}  F.-Y. Wang, \emph{Derivative formula and gradient estimates for  Gruschin type semigroups,} to appear in J. Theo. Probab.,  arXiv:1109.6738.

\bibitem{W12} F.-Y. Wang, \emph{Modified curvatures on manifolds with boundary and applications,}  arXiv:1102.3552.


\end{thebibliography}
\end{document}